\documentclass{amsart}

\usepackage{fancyhdr}
\usepackage{lastpage}
\usepackage{stmaryrd,yhmath}

\newcommand{\bdis}{\begin{displaymath}}
\newcommand{\edis}{\end{displaymath}}
\newcommand{\be}{\begin{equation}}
\newcommand{\ee}{\end{equation}}
\newcommand{\mbb}{\mathbb}
\newcommand{\mcal}{\mathcal}

\newcommand{\vp}{\varphi}

\newcommand{\zf}{\zeta\left(\frac{1}{2}+it\right)}

\theoremstyle{definition}

\newtheorem{cor}[]{Corollary}

\theoremstyle{remark}
\newtheorem{remark}[]{Remark}

\newtheorem*{mydef1}{{\bf Theorem}}

\newtheorem*{mydef7}{{\bf Question}}

\numberwithin{equation}{section}

\begin{document}

\title{Jacob's ladders, conjugate integrals, external mean-values and other properties of a multiply $\pi(T)$-autocorrelation of the
function $|\zf|^2$}

\author{Jan Moser}

\address{Department of Mathematical Analysis and Numerical Mathematics, Comenius University, Mlynska Dolina M105, 842 48 Bratislava, SLOVAKIA}

\email{jan.mozer@fmph.uniba.sk}


\begin{abstract}
In this paper we obtain a new class of transformation formulae (without an explicit presence of a derivative)
for the integrals containing products of factors $|\zf|^2$ with respect to two
components of a disconnected set on the critical line.
\end{abstract}
\maketitle

\section{Introduction}

\subsection{}
In the work of reference \cite{3} (comp. also \cite{1} and \cite{2}) we have introduced the following disconnected set
\be \label{1.1}
\Delta(n+1)=\Delta(n+1;T,U)=\bigcup_{k=0}^{n+1}[\vp_1^k(T),\vp_1^{k}(T+U)]
\ee
where
\be \label{1.2}
\begin{split}
& y=\frac 12\vp(t)=\vp_1(t);\ \vp_1^0(t)=t,\ \vp_1^1(t)=\vp_1(t),\\
& \vp_1^2(t)=\vp_1[\vp_1(t)],\dots, \vp_1^k(t)=\vp_1[\vp_1^{k-1}(t)], \dots,\ t\in [T,T+U],
\end{split}
\ee
and $\vp_1^k(t)$ stands for the $k$-th iteration of the Jacob's ladder
$$\vp_1(t),\ t\geq T_0[\vp_1].$$
The set (\ref{1.1}) has the following properties
\be \label{1.3}
\begin{split}
& t\sim \vp_1^k(t),\ \vp_1^k(T)\geq (1-\epsilon)T,\ k=0,1,\dots,n+1, \\
& \vp_1^k(T+U)-\vp_1^k(T)<\frac{1}{2n+5}\frac{T}{\ln T},\ k=1,\dots,n+1, \\
& \vp_1^k(T)-\vp_1^{k+1}(T+U)>0.18\times \frac{T}{\ln T},\ k=0,1,\dots,n, \\
& U\in\left(\left. 0,\frac{T}{\ln^2T}\right]\right. ,
\end{split}
\ee
and, in the macroscopic domain, i. e. for
\be \label{1.4}
U\in \left[ T^{1/3+\epsilon},\frac{T}{\ln^2T}\right],
\ee
we have a more detailed information about the set (\ref{1.1}), namely
\be \label{1.5}
\begin{split}
& |[\vp_1^k(T),\vp_1^{k}(T+U)]|=\vp_1^k(T+U)-\vp_1^k(T)\sim U,\ k=1,\dots,n+1, \\
& \vp_1^k(T)-\vp_1^{k+1}(T+U)\sim (1-c)\frac{T}{\ln T},\ k=0,1,\dots,n,
\end{split}
\ee
where $c$ is the Euler constant. We have that (see (\ref{1.3}))
\bdis
[\vp_1^{n+1}(T),\vp_1^{n+1}(T+U)]\prec \dots \prec [\vp_1^1(T),\vp_1^1(T+U)]\prec [T,T+U],
\edis
i. e. the segments are ordered from $[T,T+U]$ to the left.

\begin{remark}
The asymptotic behavior of the disconnected set (\ref{1.1}) is as follows: if $T\to\infty$ then the components of this set recedes unboundedly
each from other  (see (\ref{1.3}), (\ref{1.5})) and all together are receding to infinity. Hence, if $T\to\infty$ then the set (\ref{1.1}) behaves
as an one-dimensional Friedman-Hubble expanding universe.
\end{remark}

\subsection{}
Next, we have shown (see \cite{3}) that for the weighted mean-value of the integral
\be \label{1.6}
\int_T^{T+U} \prod_{k=0}^n\left|\zeta\left(\frac 12+i\vp_1^k(t)\right)\right|^2{\rm d}t,\ U\in\left(\left. 0,\frac{T}{\ln^2T}\right]\right.
\ee
the following factorization formula
\be \label{1.7}
\begin{split}
& g_{n+1}\frac{1}{U}\int_T^{T+U}\prod_{k=0}^n\left|\zeta\left(\frac 12+i\vp_1^k(t)\right)\right|^2{\rm d}t\sim \\
& \sim \prod_{l=1}^s g_l\frac{1}{U}\int_T^{T+U}\prod_{k=0}^{a_{j_l}-1}\left|\zeta\left(\frac 12+i\vp_1^k(t)\right)\right|^2{\rm d}t,\ T\to\infty
\end{split}
\ee
holds true for every fixed natural number $n$ and for every proper partition  (the partition $n+1=n+1$ is excluded)
\bdis
n+1=a_{j_1}+a_{j_2}+\dots+a_{j_s},\ a_{j_l}\in [1,n],\ l=1,\dots,s,
\edis
and
\bdis
\begin{split}
& g_l=\frac{U}{\vp_1^{a_{j_l}}(T+U)-\vp_1^{a_{j_l}}(T)},\ l=1,\dots,s, \\
& g_{n+1}=\frac{U}{\vp_1^{n+1}(T+U)-\vp_1^{n+1}(T)}.
\end{split}
\edis

\subsection{}
Next, by \cite{3}, (6.5), $n+1\to k$, we have
\bdis
t-\vp_1^k\sim k(1-c)\pi(t),\ k=0,1,\dots,n
\edis
where $\pi(t)$ is the prime-counting function. Hence
\be \label{1.8}
\begin{split}
& \frac{1}{2}+i\vp_1^k(t)=\frac{1}{2}+it-i[t-\vp_1^k(t)]\sim \\
& \sim \frac{1}{2}+it-ik(1-c)\pi(t),\ k=0,1,\dots,n.
\end{split}
\ee

\begin{remark}
By (\ref{1.8}) the arguments in the product (\ref{1.6}) performs some complicated oscillations around the sequence
\bdis
\frac{1}{2}+it-ik(1-c)\pi(t),\ k=0,1,\dots,n
\edis
of the lattice points. Based on this, the integral (\ref{1.6}) represents the multiple (for $k\geq 2$) $\pi(t)$-autocorrelation of the
function $|\zf|^2$, i. e. we have certain type of the complicated nonlinear and nonlocal interaction of the function $|\zf|^2$ with
itself.
\end{remark}

\subsection{}
After this we turn back to the formula (\ref{1.7}). This formula binds the corresponding set of integrals over the same segment $[T,T+U]$. However,
the segment $[T,T+U]$ is only one component of the disconnected set $\Delta(n+1)$ (see (\ref{1.1})). This is the reason for the following.

\begin{mydef7}
Is there some formula that binds the integral (\ref{1.6}) with the integral of the type
\bdis
\int_{\vp_1^{p(n)}(T)}^{\vp_1^{p(n)}(T+U)}\prod_{k=0}^n\left|\zeta\left(\frac 12+i\vp_1^k(u)\right)\right|^2{\rm d}u,\ 1\leq p(n)\leq n,
\edis
i. e. with the integral over the component
\bdis
[\vp_1^{p(n)}(T),\vp_1^{p(n)}(T+U)]\not= [T,T+U].
\edis
\end{mydef7}

\section{The main formula and its structure}

\subsection{}
We obtain the following theorem in the direction of our Question

\begin{mydef1}
For every disconnected set
\bdis
\Delta(2l)=\Delta(2l;T,U)=\bigcup_{k=0}^{2l} [\vp_1^k(T),\vp_1^k(T+U)],\ l=1,\dots,L_0
\edis
where $L_0\in\mbb{N}$ is an arbitrary fixed number, and for every
\bdis
U\in\left(\left. 0,\frac{T}{\ln^2T}\right]\right.
\edis
the following asymptotic transformation formula
\be \label{2.1}
\begin{split}
& \int_{\vp_1^l(T)}^{\vp_1^l(T+U)}\prod_{k=0}^{l-1}\left|\zeta\left(\frac{1}{2}+i\vp_1^k(u_l)\right)\right|^2{\rm d}u_l\sim \\
& \sim \frac{\vp_1^{2l}(T+U)-\vp_1^{2l}(T)}{\vp_1^{l}(T+U)-\vp_1^{l}(T)}\int_T^{T+U}
\prod_{k=0}^{l-1}\left|\zeta\left(\frac{1}{2}+i\vp_1^k(t)\right)\right|^2{\rm d}t,\ T\to\infty
\end{split}
\ee
holds true.
\end{mydef1}

\begin{remark}
We call the integrals that are bind by the formula (\ref{2.1}) \emph{the conjugate integrals}.
\end{remark}

Let
\bdis
\frac{1}{2}+i\gamma,\ \frac{1}{2}+i\gamma',\quad \gamma<\gamma'
\edis
be consecutive zeros of the Riemann zeta-function lying on the critical line and $l=7, T=\gamma, U=\gamma'-\gamma$. Thus, for example, the following
formula (see (\ref{2.1}))
\be \label{2.2}
\begin{split}
& \int_{\vp_1^7(\gamma)}^{\vp_1^7(\gamma')}\prod_{k=0}^6\left|\zeta\left(\frac{1}{2}+i\vp_1^k(u_7)\right)\right|^2{\rm d}u_7\sim \\
& \sim \frac{\vp_1^{14}(T+U)-\vp_1^{14}(T)}{\vp_1^{7}(T+U)-\vp_1^{7}(T)}\int_\gamma^{\gamma'}
\prod_{k=0}^6\left|\zeta\left(\frac{1}{2}+i\vp_1^k(t)\right)\right|^2{\rm d}t,\ \gamma\to\infty
\end{split}
\ee
holds true.

\begin{remark}
Nor the formula (\ref{2.2}) for seven factors and $U=\gamma'-\gamma$ is not accessible for the current methods in the theory of the Riemann
zeta-function.
\end{remark}

\subsection{}
By the continuity of the function $\vp_1^l(v)$ we have (see (\ref{2.1})) that if
\bdis
u_l=\vp_1^l(t),\ t\in [T,T+U]
\edis
then
\bdis
\vp_1^k(u_l)=\vp_1^k[\vp_1^l(t)]=\vp_1^{k+l}(t)\in [\vp_1^{k+l}(T),\vp_1^{k+l}(T+U)].
\edis
Consequently,the product
\bdis
\prod_{k=0}^{l-1}\left|\zeta\left(\frac{1}{2}+i\vp_1^k(u_l)\right)\right|^2
\edis
corresponds to the disconnected set
\be \label{2.3}
\bigcup_{k=l}^{2l-1}[\vp_1^k(T),\vp_1^k(T+U)]=\Delta(l,2l-1),
\ee
and similarly the product
\bdis
\prod_{k=0}^{l-1}\left|\zeta\left(\frac{1}{2}+i\vp_1^k(t)\right)\right|^2
\edis
corresponds to the disconnected set
\be \label{2.4}
\bigcup_{k=0}^{l-1}[\vp_1^k(T),\vp_1^k(T+U)]=\Delta(0,l-1),
\ee
where the sets (\ref{2.3}), (\ref{2.4}) are subsets of the set $\Delta(2l)$. \\

Next (comp. (\ref{1.3})), we have
\be \label{2.5}
\rho\{[\vp_1^k(T),\vp_1^k(T+U)];[\vp_1^{k+1}(T),\vp_1^{k+1}(T+U)]\}>0.17\times\pi(T)
\ee
where $\rho$ represents the distance of corresponding segments.

\begin{remark}
The formula (\ref{2.1}) controls a \emph{quasi-chaotic} behavior of the values of the function $|\zf|^2$ with respect to the disconnected
set $\Delta(2l)$ in spite of big distances separating the components of the set $\Delta(2l)$ (see (\ref{2.5})).
\end{remark}

\section{Some external mean-values}

\subsection{}
Using the mean-value theorem on the left-hand side of (\ref{2.1}) we obtain
\be \label{3.1}
\begin{split}
& \frac{1}{U}\int_T^{T+U}\prod_{k=0}^{l-1}\left|\zeta\left(\frac{1}{2}+i\vp_1^k(t)\right)\right|^2{\rm d}t\sim \\
& \sim \frac{\{\vp_1^l(T+U)-\vp_1^l(T)\}^2}{\{\vp_1^{2l}(T+U)-\vp_1^{2l}(T)\}U}\prod_{k=0}^{l-1}\left|\zeta\left(\frac{1}{2}+i\vp_1^k(\alpha_l)\right)\right|^2
\end{split}
\ee
where (see the paragraph 2.2)
\bdis
\alpha_l\in (\vp_1^l(T),\vp_1^l(T+U)),\ \alpha_l=\vp_1^l(t_l),
\edis
i. e.
\be \label{3.2}
\vp_1^k(\alpha_l)=\vp_1^{k+l}(t_l)\in (\vp_1^{k+l}(T),\vp_1^{k+l}(T+U)).
\ee
Hence, by (\ref{3.1}) and (\ref{3.2}) we have the following
\begin{cor}
There are the values
\bdis
\tau_k=\tau_k(T,U,l)\in (\vp_1^k(T),\vp_1^k(T+U)),\ k=l,\dots,2l-1
\edis
such that
\be \label{3.3}
\begin{split}
& \frac{1}{U}\int_T^{T+U}\prod_{k=0}^{l-1}\left|\zeta\left(\frac{1}{2}+i\vp_1^k(t)\right)\right|^2{\rm d}t\sim \\
& \sim \frac{\{\vp_1^l(T+U)-\vp_1^l(T)\}^2}{\{\vp_1^{2l}(T+U)-\vp_1^{2l}(T)\}U}\prod_{k=l}^{2l-1}\left|\zeta\left(\frac{1}{2}+i\tau_k\right)\right|^2
\end{split}
\ee
where
\bdis
U\in \left(\left. 0,\frac{T}{\ln^2T}\right]\right.,\ l=1,\dots,L_0,\ T\to\infty.
\edis
\end{cor}

\begin{remark}
Since:
\begin{itemize}
\item[(a)] the integral
\bdis
\int_T^{T+U}\prod_{k=0}^{l-1}\left|\zeta\left(\frac{1}{2}+i\vp_1^k(t)\right)\right|^2{\rm d}t
\edis
corresponds to the disconnected set $\Delta(0,l-1)$, (see (\ref{2.4})),
\item[(b)] the product
\bdis
\prod_{k=l}^{2l-1}\left|\zeta\left(\frac{1}{2}+i\tau_k\right)\right|
\edis
corresponds to the disconnected set $\Delta(l,2l-1)$, (see (\ref{2.3})),
\item[(c)] the sets $\Delta(0,l-1)$ and $\Delta(l,2l-1)$ are separated by the big distance
\bdis
\rho\{\Delta(0,l-1);\Delta(l,2l-1)\}>0.17\times \pi(T)
\edis
(see (\ref{2.3}), (\ref{2.4})),
\end{itemize}
it is quite natural to call the right-hand side of the equation (\ref{3.3}) \emph{the external mean-value} of the integral on the left-hand
side.
\end{remark}

\subsection{}
Next, by the similar way, we obtain the following
\begin{cor}
There are the values
\bdis
\tau_k=\tau_k(T,U,l)\in (\vp_1^k(T),\vp_1^k(T+U)),\ k=0,1,\dots,l-1
\edis
such that
\be \label{3.4}
\begin{split}
& \frac{1}{\vp_1^l(T+U)-\vp_1^l(T)}\int_{\vp_1^l(T)}^{\vp_1^l(T+U)}\prod_{k=0}^{l-1}
\left|\zeta\left(\frac{1}{2}+i\vp_1^k(u_l)\right)\right|^2{\rm d}u_l\sim \\
& \sim \frac{\{\vp_1^{2l}(T+U)-\vp_1^{2l}(T)\}U}{\{\vp_1^l(T+U)-\vp_1^l(T)\}^2}\prod_{k=0}^{l-1}\left|\zeta\left(\frac{1}{2}+i\tau_k\right)\right|^2,
\end{split}
\ee
where
\bdis
U\in \left(\left. 0,\frac{T}{\ln^2T}\right]\right.,\ l=1,\dots,L_0,\ T\to\infty.
\edis
\end{cor}

\begin{remark}
The formula (\ref{3.4}) gives us the second variant of the external mean-value theorem.
\end{remark}

\section{Other properties of the distribution of the values of $|\zf|$ with respect to the disconnected set $\Delta(2l)$}

\subsection{}
Similarly to (\ref{3.3}), (\ref{3.4}), we obtain the following formula
\be \label{4.1}
\begin{split}
& \prod_{k=0}^l\left|\zeta\left(\frac 12+i\tau_k\right)\right|\sim \\
& \sim \frac{\vp_1^l(T+U)-\vp_1^l(T)}{\sqrt{\{\vp_1^{2l}(T+U)-\vp_1^{2l}(T)\}U}}\prod_{k=l}^{2l-1}\left|\zeta\left(\frac 12+i\tau_k\right)\right|,
\end{split}
\ee
where
\bdis
\tau_k\in (\vp_1^k(T),\vp_1^k(T+U)),\ k=0,1,\dots,2l-1.
\edis
Next, we obtain from (\ref{4.1}) the following
\begin{cor}
\be \label{4.2}
\begin{split}
& G_0^{l-1}\left[\left|\zeta\left(\frac 12+i\tau_k\right)\right|\right]\sim \\
& \sim \left\{\frac{\vp_1^l(T+U)-\vp_1^l(T)}{\sqrt{\{\vp_1^{2l}(T+U)-\vp_1^{2l}(T)\}U}}\right\}^{1/l}
G_l^{2l-1}\left[\left|\zeta\left(\frac 12+i\tau_k\right)\right|\right],\ T\to\infty
\end{split}
\ee
where the following symbols
\be \label{4.3}
\begin{split}
& G_0^{l-1}\left[\left|\zeta\left(\frac 12+i\tau_k\right)\right|\right]=
\left\{\prod_{k=0}^{l-1}\left|\zeta\left(\frac 12+i\tau_k\right)\right| \right\}^{1/l}, \\
& G_l^{2l-1}\left[\left|\zeta\left(\frac 12+i\tau_k\right)\right|\right]=
\left\{\prod_{k=l}^{2l-1}\left|\zeta\left(\frac 12+i\tau_k\right)\right| \right\}^{1/l}
\end{split}
\ee
stand for the geometric means.
\end{cor}

\subsection{}
Since (see (\ref{4.3}))
\be \label{4.4}
\frac{G_0^{l-1}}{G_l^{2l-1}}=\bar{G}_0^{l-1}
\left[\frac{\left|\zeta\left(\frac{1}{2}+i\tau_k\right)\right|}{\left|\zeta\left(\frac{1}{2}+i\tau_{k+l}\right)\right|}\right],
\ee
and we have for arithmetic and geometric means (for example)
\be \label{4.5}
\bar{x}_A\geq \bar{x}_G;\ \bar{x}_A=\frac{1}{n}\sum_{i=1}^n x_i,\ \bar{x}_G=\sqrt[n]{\prod_{i=1}^n x_i},\ x_i>0.
\ee
Then we obtain from (\ref{4.2})-(\ref{4.4}) the formula
\bdis
\bar{G}_0^{l-1}
\left[\frac{\left|\zeta\left(\frac{1}{2}+i\tau_k\right)\right|}{\left|\zeta\left(\frac{1}{2}+i\tau_{k+l}\right)\right|}\right]\sim
\left\{\frac{\vp_1^l(T+U)-\vp_1^l(T)}{\sqrt{\{\vp_1^{2l}(T+U)-\vp_1^{2l}(T)\}U}}\right\}^{1/l}=\Omega_l.
\edis
Next, from the inequality
\bdis
\bar{G}^{l-1}_0>(1-\epsilon)\Omega_l,\ T\to\infty
\edis
we obtain that (see (\ref{4.5}))
\be \label{4.6}
\frac{1}{l}\sum_{k=0}^{l-1}\frac{\left|\zeta\left(\frac{1}{2}+i\tau_k\right)\right|}{\left|\zeta\left(\frac{1}{2}+i\tau_{k+l}\right)\right|}
>(1-\epsilon)\Omega_l.
\ee
The numbers $(\tau_0,\tau_1,\dots,\tau_{l-1})$ may be ordered by $l!$-ways in the product
\bdis
\prod_{k=0}^{l-1}\frac{\left|\zeta\left(\frac{1}{2}+i\tau_k\right)\right|}{\left|\zeta\left(\frac{1}{2}+i\tau_{k+l}\right)\right|},
\edis
and the same holds for the sequence of numbers $(\tau_l,\dots,\tau_{2l-1})$. Therefore we have $(l!)^2$ inequalities of the type (\ref{4.6}). In this sense
we use the symbol
\bdis
\left\{\sum_{k=0}^{l-1}\frac{\left|\zeta\left(\frac{1}{2}+i\tau_k\right)\right|}{\left|\zeta\left(\frac{1}{2}+i\tau_{k+l}\right)\right|}\right\}_m,\ m=1,\dots,(l!)^2.
\edis
Hence, we obtain from (\ref{4.6}) the following
\begin{cor}
We have $(l!)^2$ inequalities
\bdis
\frac{1}{l}
\left\{\sum_{k=0}^{l-1}\frac{\left|\zeta\left(\frac{1}{2}+i\tau_k\right)\right|}{\left|\zeta\left(\frac{1}{2}+i\tau_{k+l}\right)\right|}\right\}_m>
(1-\epsilon)
\left\{\frac{\vp_1^l(T+U)-\vp_1^l(T)}{\sqrt{\{\vp_1^{2l}(T+U)-\vp_1^{2l}(T)\}U}}\right\}^{1/l},
\edis
for $\tau_0,\tau_1,\dots,\tau_{2l-1}$, where
\bdis
m=1,\dots,(l!)^2,\ l=1,\dots,L_0,\ U\in \left(\left. 0,\frac{T}{\ln^2T}\right]\right.,\ l=1,\dots,L_0,\ T\to\infty.
\edis
\end{cor}

\begin{remark}
There are certain multiplicative effects also in the genetics, among the polygenic systems, and consequently the geometric means is used there, see,
for example, \cite{4}, pp. 336, 337.
We also note that we have used the formula for multiplication of independent variables as a motivation for our paper \cite{3} .
\end{remark}

\section{Remarks about essential influence of the Riemann hypothesis on the sequence $\{\vp_1^k(T+U)-\vp_1^k(T)\}_{k=1}^{L_0}$}

\subsection{}
Let us remind that in the macroscopic case (\ref{1.4}) we have the asymptotic formula (see (\ref{1.5}))
\be \label{5.1}
\vp_1^k(T+U)-\vp_1^k(T)\sim U,\ k=1,\dots,L_0.
\ee
In connection with (\ref{5.1}) we ask the question: what is the influence of the Riemann hypothesis on measures of the segments
\bdis
[\vp_1(T),\vp_1(T+U)]
\edis
in the case (comp. (\ref{1.4}))
\be \label{5.2}
U\in (0,T^{1/3-\epsilon_0}],
\ee
for example, in the case $\epsilon_0=\frac{1}{12}$, i. e.
\bdis
U\in (0,T^{1/4}].
\edis

First of all we have, on the Riemann hypothesis, that (see \cite{5}, p. 300)
\be \label{5.3}
\zf=\mcal{O}\left( t^{\frac{A}{\ln\ln t}}\right),\ t\to\infty,
\ee
i. e.
\be \label{5.4}
\zf=\mcal{O}\left( T^{\frac{A}{\ln\ln T}}\right),\ t\in [(1-\epsilon)T,T+U]
\ee
(comp. (\ref{1.3}) and \cite{3}, (6.17)). Next we obtain for (\ref{5.2}) from our formula (see \cite{2}, (2.5))
\bdis
\begin{split}
& \int_T^{T+V}\left|\zf\right|^2{\rm d}t\sim [\vp_1(T+V)-\vp_1(T)]\ln T, \\
& V\in \left(\left. 0,\frac{T}{\ln T}\right.\right],
\end{split}
\edis
by (\ref{5.4}) that
\be \label{5.5}
\begin{split}
& \vp_1^1(T+U)-\vp_1^1(T)=\mcal{O}\left(\frac{U}{\ln T}T^{\frac{2A}{\ln\ln T}}\right), \\
& \vp_1^2(T+U)-\vp_1^2(T)=\mcal{O}\left(\frac{U}{\ln^2 T}T^{2\frac{2A}{\ln\ln T}}\right), \\
& \vdots \\
& \vp_1^{L_0}(T+U)-\vp_1^{L_0}(T)=\mcal{O}\left(\frac{U}{\ln^{L_0} T}T^{L_0\frac{2A}{\ln\ln T}}\right).
\end{split}
\ee
Since
\be \label{5.6}
T^{L_0\frac{2A}{\ln\ln T}}=T^{\frac{2L_0A}{\sqrt{\ln\ln T}}\frac{1}{\sqrt{\ln\ln T}}}<T^{\frac{1}{\sqrt{\ln\ln T}}},
\ee
then by (\ref{5.5}), (\ref{5.6}) we obtain the following

\begin{remark}
On the Riemann hypothesis the following estimates hold true
\be \label{5.7}
\begin{split}
& U\in (0,T^{1/3-\epsilon}] \ \Rightarrow \ \vp_1^k(T+U)-\vp_1^k(T)=
\mcal{O}\left( UT^{\frac{1}{\sqrt{\ln\ln T}}}\right),\\
& k=1,\dots,L_0.
\end{split}
\ee
For example, if $U=1$ then on Riemann hypothesis we have that
\bdis
\vp_1^k(T+1)-\vp_1^k(T)=\mcal{O}\left( UT^{\frac{1}{\sqrt{\ln\ln T}}}\right),\ k=1,\dots,L_0
\edis
either for
\bdis
L_0=S=10^{10^{10^{34}}}
\edis
($S$ is the Skeewes' constant).
\end{remark}

\subsection{}
In the general case (with or without the Riemann hypothesis) we have (comp. (\ref{5.3}), (\ref{5.4}))
\bdis
\zf=\mcal{O}(t^{1/6-\epsilon})=\mcal{O}(T^{1/6-\epsilon}),\ t\in [(1-\epsilon)T,T+U],\ T\to\infty,
\edis
and consequently we obtain (comp. (\ref{5.5}))
\bdis
\begin{split}
& \vp_1^1(T+1)-\vp_1^1(T)=\mcal{O}(T^{2(1/6-\epsilon)})=\mcal{O}(T^{1/3-2\epsilon}), \\
& \vp_1^2(T+1)-\vp_1^2(T)=\mcal{O}(T^{4(1/6-\epsilon)})=\mcal{O}(T^{2/3-4\epsilon}).
\end{split}
\edis

\begin{remark}
In the general case we are able to guarantee only that
\be \label{5.8}
\vp_1^1(T+1)-\vp_1^1(T)\in (0,T^{1/3-\epsilon_0}],\ \epsilon\leq \frac{\epsilon_0}{2}.
\ee
Hence, the comparison of (\ref{5.7}), $U=1$, with (\ref{5.8}) shows the essential influence of the Riemann hypothesis on our subject.
\end{remark}

\section{The proof of Theorem}

\subsection{}
By using our formula (see \cite{2}, (9.1))
\bdis
\tilde{Z}^2(t)=\frac{{\rm d}\vp_1(t)}{{\rm d}t}
\edis
we obtain (see (\ref{1.2}))
\bdis
\begin{split}
& \int_T^{T+U}\prod_{k=0}^n \tilde{Z}^2[\vp_1^k(t)]{\rm d}t=\\
& = \int_T^{T+U}\tilde{Z}^2[\vp_1^n(t)]\tilde{Z}^2[\vp_1^{n-1}(t)]\cdots\tilde{Z}^2[\vp_1^1(t)]\tilde{Z}^2[t]{\rm d}t=\\
& =\int_T^{T+U}\tilde{Z}^2[\vp_1^{n-1}(\vp_1^1(t))]\tilde{Z}^2[\vp_1^{n-2}(\vp_1^1(t))]\cdots\tilde{Z}^2[\vp_1^1(t)]\frac{{\rm d}\vp_1^1(t)}{{\rm d}t}{\rm d}t=\\
& =\int_{\vp_1^1(T)}^{\vp_1^1(T+U)}\tilde{Z}^2[\vp_1^{n-1}(u_1)]\tilde{Z}^2[\vp_1^{n-2}(u_1)]\cdots\tilde{Z}^2[\vp_1^1(u_1)]\tilde{Z}^2[u_1]{\rm d}u_1=\\
& =\int_{\vp_1^1(T)}^{\vp_1^1(T+U)}\tilde{Z}^2[\vp_1^{n-2}(\vp_1^1(u_1))]\cdots\tilde{Z}^2[\vp_1^1(u_1)]\frac{{\rm d}\vp_1^1(u_1)}{{\rm d}u_1}{\rm d}u_1= \\
& = \int_{\vp_1^2(T)}^{\vp_1^2(T+U)}\tilde{Z}^2[\vp_1^{n-2}(u_2)]\cdots\tilde{Z}^2[u_2]{\rm d}u_2=\dots=\\
& =\int_{\vp_1^l(T)}^{\vp_1^l(T+U)}\tilde{Z}^2[\vp_1^{n-l}(u_l)]\cdots\tilde{Z}^2[\vp_1^0(u_l)]{\rm d}u_l,\ l=1,\dots,n,
\end{split}
\edis
i. e. the following formula
\be \label{6.1}
\begin{split}
& \int_T^{T+U}\prod_{k=0}^n\tilde{Z}^2[\vp_1^k(t)]{\rm d}t=\int_{\vp_1^l(T)}^{\vp_1^l(T+U)}\prod_{k=0}^{n-l}\tilde{Z}^2[\vp_1^k(u_l)]{\rm d}u_l, \\
& l=1,\dots,n
\end{split}
\ee
holds true.

\subsection{}
Let us remind that (see \cite{3}, (6.14))
\be \label{6.2}
\begin{split}
& \tilde{Z}^2(t)=\frac{Z^2(t)}{2\Phi'_\vp[\vp(t)]}=\frac{\left|\zf\right|^2}{\left\{ 1+\mcal{O}\left(\frac{\ln\ln T}{\ln T}\right)\right\}\ln t}, \\
& t\in [T,T+U],\ U\in\left(\left. 0,\frac{T}{\ln T}\right]\right., \\
& (\vp_1^l(T),\vp_1^l(T+U))\subset (\vp_1^{n+1}(T),T+U).
\end{split}
\ee
Putting (\ref{6.2}) into (\ref{6.1}) and using the mean-value theorem on both integrals in (\ref{6.1}) we obtain the following formula (comp. \cite{3}, (6.17))
\be \label{6.3}
\begin{split}
& \int_T^{T+U}\prod_{k=0}^n \left|\zeta\left(\frac{1}{2}+i\vp_1^k(t)\right)\right|^2{\rm d}t\sim \\
& \sim \ln^lT\int_{\vp_1^l(T)}^{\vp_1^l(T+U)}\prod_{k=0}^{n-l}\left|\zeta\left(\frac{1}{2}+i\vp_1^k(u)\right)\right|^2{\rm d}u,\ l=1,\dots,n,\ T\to\infty.
\end{split}
\ee
Next, the formula (see \cite{3}, (3.1))
\be \label{6.4}
\begin{split}
& \int_T^{T+U}\prod_{k=0}^n\left|\zeta\left(\frac{1}{2}+i\vp_1^k(t)\right)\right|^2{\rm d}t\sim
\{ \vp_1^{n+1}(T+U)-\vp_1^{n+1}(T)\}\ln^{n+1}T; \\
& \ln^{n+1}T=\ln^{(l-1)+1}T\ln^{(n-l)+1}T
\end{split}
\ee
together with the formula (\ref{6.3}) gives the following asymptotic equality
\bdis
\frac{\int_T^{T+U}}{\vp_1^{n+1}(T+U)-\vp_1^{n+1}(T)}\sim \frac{\int_T^{T+U}}{\int_{\vp_1^l(T)}^{\vp_1^l(T+U)}}\frac{\int_T^{T+U}}{\int_{\vp_1^{n+1-l}(T)}^{\vp_1^{n+1-l}(T+U)}},
\edis
i. e.
\be \label{6.5}
\begin{split}
& \{\vp_1^{n+1}(T+U)-\vp_1^{n+1}(T)\}\int_T^{T+U}\prod_{k=0}^n\left|\zeta\left(\frac{1}{2}+i\vp_1^k(t)\right)\right|^2{\rm d}t\sim \\
& \sim \int_{\vp_1^l(T)}^{\vp_1^l(T+U)}\prod_{k=0}^{n-l}\left|\zeta\left(\frac{1}{2}+i\vp_1^k(u)\right)\right|^2{\rm d}u \times \\
& \times
\int_{\vp_1^{n+1-l}(T)}^{\vp_1^{n+1-l}(T+U)}\prod_{k=0}^{l-1}\left|\zeta\left(\frac{1}{2}+i\vp_1^k(v)\right)\right|^2{\rm d}v,\
l=1,\dots,n,\ T\to\infty.
\end{split}
\ee

\subsection{}
Next, in the case
\bdis
n-l=l-1 \ \Rightarrow \ n=2l-1,
\edis
we obtain that (see (\ref{6.4}), (\ref{6.5}))
\bdis
\begin{split}
& \left\{\int_{\vp_1^l(T)}^{\vp_1^l(T+U)}\prod_{k=0}^{l-1}\left|\zeta\left(\frac{1}{2}+i\vp_1^k(u_l)\right)\right|^2{\rm d}u_l\right\}^2\sim \\
& \sim \{\vp_1^{2l}(T+U)-\vp_1^{2l}(T)\}\int_T^{T+U}\prod_{k=0}^{2l-1}\left|\zeta\left(\frac{1}{2}+i\vp_1^k(t)\right)\right|^2{\rm d}t\sim \\
& \sim \{\vp_1^{2l}(T+U)-\vp_1^{2l}(T)\}^2\ln^{2l}T,
\end{split}
\edis
i. e. the following formula holds true
\be \label{6.6}
\begin{split}
& \int_{\vp_1^l(T)}^{\vp_1^l(T+U)}\prod_{k=0}^{l-1}\left|\zeta\left(\frac{1}{2}+i\vp_1^k(u_l)\right)\right|^2{\rm d}u_l\sim \\
& \sim \{\vp_1^{2l}(T+U)-\vp_1^{2l}(T)\}\ln^{l}T.
\end{split}
\ee
Consequently, we obtain from (\ref{6.6}) by (\ref{6.4}), in the case $n=l-1$, the formula
\bdis
\begin{split}
& \int_{\vp_1^l(T)}^{\vp_1^l(T+U)}\prod_{k=0}^{l-1}\left|\zeta\left(\frac{1}{2}+i\vp_1^k(u_l)\right)\right|^2{\rm d}u_l\sim \\
& \sim \frac{\vp_1^{2l}(T+U)-\vp_1^{2l}(T)}{\vp_1^{l}(T+U)-\vp_1^{l}(T)}
\int_T^{T+U}\prod_{k=0}^{l-1}\left|\zeta\left(\frac{1}{2}+i\vp_1^k(t)\right)\right|^2{\rm d}t
\end{split}
\edis
that verifies (\ref{2.1}).

\thanks{I would like to thank Michal Demetrian for his help with electronic version of this paper.}

\end{document}